\documentclass[12pt,a4paper]{article}
\def\ba{\begin{array}}
\def\ea{\end{array}}
\def\be{\begin{equation}}
\def\ee{\end{equation}}
\def\lbl{\label}
\def \rf  {(\ref}

\def\x0{\x_0}
\def\x1{\x_1}

\def\cd{\cal {D}}
\def\cg{{\cal {G}}}
\def\tcg{{\tilde{\cal G}}}
\def\abg{  \alpha^2 + \beta^2 -  \gamma^2}

\begin{document}
\author{L. Hlavat\'y, L. \v Snobl
\thanks{{Email: hlavaty@br.fjfi.cvut.cz}, {snobl@newton.fjfi.cvut.cz}} 
\\ {\it Faculty of Nuclear Sciences and Physical Engineering,} 
\\{ \it Czech Technical University,} 
\\ {\it B\v rehov\'a 7, 115 19 Prague 1, Czech Republic}}

\title{Classification of 6--dimensional real Manin triples }
\date{March 27, 2002}

\maketitle
\bibliographystyle{unsrt}

\abstract{We present a complete list of 6--dimensional Manin triples 
or, equivalently, of 3--dimensional Lie bialgebras. We start from the well known classification of 3--dimensional real Lie algebras and assume the canonical bilinear form on the 6-dimensional Drinfeld double. Then we solve the Jacobi identities for the dual algebras. Finally we find mutually non--isomorphic Manin triples.
The complete list consists of 78 classes of  Manin triples, or 44 Lie bialgebras if one considers dual bialgebras equivalent. 
}
\newpage\section{Introduction}
In recent years, the study of T--duality in string theory has led to 
discovery of Poisson--Lie T--dual sigma models. Klim\v{c}\'{\i}k and \v{S}evera 
have found a procedure allowing to construct such models 
from a given Manin triple  
$({\cal D},{\cal G},{\tilde{\cal G}})$ , 
i.e. a decomposition of a Lie algebra 
${\cal D}$ 
into two maximally isotropic subalgebras ${\cal G},{\tilde{\cal G}}$.
The construction of the Poisson--Lie T--dual sigma models is described in \cite{klse:dna} and \cite{kli:pltd}. 
The models have target spaces in the Lie groups $G$ and $\tilde G$ and are defined by the Lagrangians
\be {\cal L}= E_{ij}(g)(g^{-1}\partial_- g)^i(g^{-1}\partial_+ g)^j \ee
\be \tilde{\cal L}=\tilde E_{ij}(\tilde g)(\tilde g^{-1}\partial_- \tilde g)^i(\tilde g^{-1}\partial_+ \tilde g)^j \ee
where the matrices $E(g)$ and $\tilde E(\tilde g)$ are constructed from 
a constant invertible matrix $E(e)$ by virtue of 
the adjoint representation of the group $G$ resp $\tilde G$ on $\cd$.
It implies that any pair of Poisson--Lie T--dual sigma models is given (up to the constant matrix $E$) by the corresponding Manin triple 
and that's why it is interesting and useful to classify these structures. 

One can easily see that the dimension of the Lie algebra ${\cal D}$ 
must be even. In the dimension two ${\cal D}$ must be abelian and there is just one Manin triple $({\cal D},{\cal G},{\tilde{\cal G}})
\equiv ({\cal D},{\tilde{\cal G}},{\cal G})$ where $dim{{\cal G}}=dim{\tilde{\cal G}}=1$. The classification of Manin triples for the 
four--dimensional Lie algebras together with the pairs of dual models was given in \cite{hlasno:pltdm2dt}. In this paper we are 
going to classify the Manin triples of the 6--dimensional real Lie algebras.

Important steps in this direction were made in \cite{iranci} where a 
list of possible maximally isotropic subalgebras of the 6--dimensional  
Lie algebras can be found. 
It turns out that the subalgebras don't specify the Manin triple completely. 
For certain algebras there exist several rather different possible pairings, 
allowing to construct different Manin triples. In the present paper, we present a complete list of real
6--dimensional Manin triples, i.e. we give not only 
the possible subalgebras, but also the corresponding ad--invariant form (i.e 
we write dual bases of the algebras with respect to this form and their Lie brackets). The complex solvable Manin triples were classified in \cite{fof:sos}.

As we shall see Manin triples are equivalent to Lie bialgebras and the classification of the three--dimensional Lie bialgebras 
(i.e. six--dimensional Manin triples) was given in \cite{gom:ctd}. Our classification was done independently without knowledge of \cite{gom:ctd}. The consequent comparison proved that the results are identical even though 
we have started from a different description of the three--dimensional algebras and used a completely different method. It means that the present work can be considered as 
an independent check of  \cite {gom:ctd} with the results expressed
in a different form, namely as Manin triples. 

In the following sections, we firstly recall the definitions 
of Manin triple, Drinfeld double and Lie bialgebra, then briefly explain the 
approach we have used to find all algebras of 6--dimensional Drinfeld doubles,
and finally give a complete list of 
all 6--dimensional 
Manin triples.

\section{Manin triples, Drinfeld doubles, Lie bialgebras}
The Drinfeld double $D$ is defined as a Lie group such that its Lie algebra 
$\cd$ equipped by a symmetric ad--invariant nondegenerate bilinear form 
$\langle .,.\rangle $ can be decomposed into a pair of maximally isotropic 
subalgebras $\cg$, $\tcg$ such that $\cd$ as a vector space is the direct 
sum of $\cg$ and $\tcg$. This ordered triple of algebras $(\cd,\cg$,$\tcg)$ is called Manin triple. 

One can see that the dimensions of the subalgebras are equal and that bases 
$\{X_i\}, \{\tilde X_i\}$ in the subalgebras can be chosen so that
\be \langle X_i,X_j\rangle =0,\  \langle X_i,\tilde X^j\rangle =\langle 
\tilde X^j,X_i\rangle =\delta_i^j,\  \langle \tilde X^i,\tilde X^j\rangle =0.\lbl{brackets}\ee
This canonical form of the bracket is invariant with respect to the transformations 
\be X_i'=X_k A^k_i,\ \tilde X^{'j}=(A^{-1})^j_k \tilde X^k. \lbl{tfnb}\ee
Due to the ad-invariance of $\langle .,.\rangle $ the algebraic structure of $\cd$ is
\[ [X_i,X_j]={f_{ij}}^k X_k,\ [\tilde X^i,\tilde X^j]={\tilde {f^{ij}}_k} \tilde X^k,\]
\be [X_i,\tilde X^j]={f_{ki}}^j \tilde X^k +{\tilde {f^{jk}}_i} X_k. \lbl{liebd}\ee

It is clear that to any Manin triple $({\cal D},{\cal G},{\tilde{\cal G}})$   one can construct the dual one by 
interchanging $\cg \leftrightarrow \tcg$, i.e. interchanging the structure constants 
$ {f_{ij}}^k \leftrightarrow {\tilde {f^{ij}}_k}$. All properties of Lie algebras 
(the nontrivial being the Jacobi identities) remain to be satisfied. On the other hand for given Drinfeld double more than two 
Manin triples can exist.

One can rewrite the structure 
of a  Manin triple also in another,
equivalent, but for certain considerations more suitable, form of Lie bialgebra. 

A Lie bialgebra is a Lie algebra $g$ equipped also by a Lie cobracket\footnote{Summation index is suppressed} 
 $\delta:g \rightarrow g \otimes g: \delta(x)=  \sum x_{[1]} \otimes x_{[2]}$ such that 
 \begin{eqnarray}
\sum x_{[1]} \otimes x_{[2]} & = & - \sum x_{[2]} \otimes x_{[1]}, \\
 (id \otimes \delta) \circ \delta (x) & + & {\rm cyclic \, permutations \, of \, tensor \, indices} = 0, \label{dji1} \\
\nonumber \delta([x,y]) & = & \sum [x,y_{[1]}] \otimes y_{[2]} + y_{[1]} \otimes [x,y_{[2]}] - \\
 & - & [y,x_{[1]}] \otimes x_{[2]} - x_{[1]} \otimes [ y, x_{[2]} ] \label{mji1}
\end{eqnarray}
(for detailed account on Lie bialgebras see e.g. \cite{Drinfeld} or \cite{Majid}, Chapter 8).

The correspondence between a Manin triple and a Lie bialgebra can now be formulated in the following way.
Because both subalgebras $\cg$, $\tcg$ of the Manin triple are of the same dimension
and are connected by nondegenerate pairing, it is natural to consider $\tcg$ as a dual $\cg^{*}$ to 
$\cg$ and to use the Lie bracket in $\tcg$ to define the Lie cobracket in $\cg$; 
$\delta(x) $ is given by $\langle \delta(x), \tilde{y} \otimes \tilde{z} \rangle = 
\langle x, [\tilde{y},\tilde{z}] \rangle,$ $\forall \tilde{y},\tilde{z} \in \cg^{*}$,
 i.e. $\delta(X_{i})= \tilde {f^{jk}_{i}}  X_{j} \otimes X_{k} $.
The Jacobi identities in $\tcg$ 
\be\label{jidd}
\tilde {f^{kl}_{m}}  \tilde {f^{ij}_{l}} + \tilde {f^{il}_{m}}  \tilde {f^{jk}_{l}} + \tilde {f^{jl}_{m}}  \tilde {f^{ki}_{l}} =0
\ee 
are then equivalent to the property of cobracket (\ref{dji1}) and 
the $\tcg$--component of the mixed Jacobi identities \footnote{The Jacobi identities $[X_i,[\tilde X^j,\tilde X^k]]+ 
{\rm cyclic} = 0$ lead to both (\ref{mjid})  (terms proportional to $\tilde X^l$) and (\ref{jidd})  (terms
proportional to $X_{l}$).} 
\be \label{mjid}
 \tilde {f^{jk}}_{l} {f_{mi}}^{l} +   \tilde {f^{kl}}_{m} {f_{li}}^{j} + 
  \tilde {f^{jl}}_{i} {f_{lm}}^{k} +   \tilde {f^{jl}}_{m} {f_{il}}^{k} + 
  \tilde {f^{kl}}_{i} {f_{lm}}^{j}   =0
\ee
are equivalent to (\ref{mji1}).  

From now on, we will use the formulation in terms of Manin triples, 
Lie bialgebra formulation of all results can be easily 
derived from it. We also consider 
only algebraic structure, the Drinfeld doubles as the Lie groups can
be obtained in principle by means of exponential map and usual
theorems about relation between Lie groups and Lie algebras 
apply, e.g. there is a one to one correspondence between (finite--dimensional) Lie 
algebras and connected and simply connected Lie groups. The group
structure of the Drinfeld double can be deduced e.g. by taking 
matrix exponential of adjoint representation of its algebra.

\section{Method of classification}

In this section we present the approach we have used to find all
6--dimensional Manin triples, i.e. 3--dimensional Lie
bialgebras.

Starting point for our computations is the well known classification 
of 3--dimensional 
real Lie algebras (see e.g. \cite{Landau} or \cite{iranci}). 
Non--isomorphic Lie algebras are written 
in 11 classes, traditionally known as Bianchi algebras. 
Their commutation relations are:
\be 
[X_1,X_2]=-a X_2+n_3 X_3, \; [X_2,X_3]=n_1 X_1, \; [X_3,X_1]=n_2 X_2 + a X_3,
\label{bian}\ee
where the parameters $a,n_1,n_2,n_3$ have the following values
\medskip

\begin{center}
\begin{tabular}{|c|r|r|r|r|}
\hline
{\bf Class} & $a$ & $n_1$ & $n_2$ & $n_3$ \\
\hline
$I$ & 0 & 0 & 0 & 0 \\
$II$ & 0 & 1 & 0 & 0 \\
$VII_0$ & 0 & 1 & 1 & 0 \\
$VI_0$ & 0 & 1 & -1 & 0 \\
$IX$ & 0 & 1 & 1 & 1 \\
$VIII$ & 0 & 1 & 1 & -1 \\
$V$ & 1 & 0 & 0 & 0 \\
$IV$ & 1 & 0 & 0 & 1 \\
$VII_a   \, (a>0)$ & $a$ & 0 & 1 & 1 \\
$III$ & 1 & 0 & 1 & -1 \\
$VI_a \, (a>0,a \neq 1)$ & $a$ & 0 & 1 & -1 \\   
\hline
\end{tabular}
\end{center}

Therefore the 1st subalgebra $\cg$  
of the Manin triple $\cd$ must be one of the Bianchi algebras given above and we can choose its basis so that the Lie brackets are of 
the form \rf{bian}).
In the 2nd subalgebra $\tcg$ we choose the dual basis $\tilde{X}^{i}$ so that \rf{brackets}) holds,  and treat nine  
independent components of structure coefficients ${\tilde{f}^{ij}}_{k}$ of the 2nd subalgebra $\tcg $
in the basis $\tilde{X}^{i}$ as unknowns. We cannot assume that the ${\tilde{f}^{ij}}_{k}$ are of the form \rf{bian}) as well because 
it can be incompatible with \rf{brackets}). Then we solve the mixed Jacobi identities 
(\ref{mjid}) 
(these relations form a system of linear equations in  ${\tilde{f}^{ij}}_{k}$) 
and the Jacobi identities for the dual algebra (\ref{jidd}) 
(i.e. quadratic in ${\tilde{f}^{ij}}_{k}$). 

As a result, we have found all  structure coefficients of $\tcg$
consistent with the definition of Manin triple and
the next step was to determine the Bianchi classes of obtained algebras $\tcg$. Finally we have found the 
the non--isomorphic Manin triples
by considering Manin triples connected by the transformations (\ref{tfnb}) 
(i.e. change of basis in $\cg$ accompanied by the dual change 
of basis in $\tcg$ with respect to $\langle,\rangle$) as equivalent and 
choosing one representant in each equivalence class.

In computations computer algebra systems Maple V and Mathematica 4 were 
independently used for manipulating expressions and solving sets of 
linear and quadratic equations, their results were checked one against the other.

Before listing our results, we shall give an example showing the
progress of computation in some detail.
\smallskip

{\noindent \bf Example:} Let us consider the algebra $VIII$, i.e. $\cg=sl(2,\Re)$.
$$ [X_1,X_2]=- X_3, \, [X_2,X_3]=X_1, \, [X_3,X_1]= X_2. $$ 
When one solves the mixed Jacobi identities (\ref{mjid}), 
he finds that the 2nd subalgebra must have the form
$$ [ \tilde{X}^1, \tilde{X}^2 ]= -\alpha \tilde{X}^1 + \beta  \tilde{X}^2, 
\, [ \tilde{X}^2, \tilde{X}^3]= \gamma \tilde{X}^2 + \alpha \tilde{X}^3 , 
\, [ \tilde{X}^3, \tilde{X}^2]= -\gamma \tilde{X}^1 - \beta \tilde{X}^3.$$ 
The Jacobi identities in the 2nd subalgebra (\ref{jidd}) in this case don't
impose any further condition on the structure constants $\tilde {f^{ij}}_{k}$, 
i.e. we have already found the structure of all possible 2nd subalgebras $\tcg$ 
in the Manin triple.

Next we find the Bianchi forms of $\tcg$. 
It turns out that the 2nd algebra is of the Bianchi type $I$ ($\tilde {f^{ij}_{k}} = 0 $) if 
$\alpha=\beta=\gamma=0$ and of type $V$ otherwise.

Then we find values of $\tilde {f^{ij}}_{k}$ that allow transformation (\ref{tfnb}) 
leading to the rescaled Bianchi form $V$ of the 2nd subalgebra $\tcg$ and leaving the Bianchi form of
the 1st subalgebra $sl(2,\Re)$ invariant. This is possible only for 
$$ \abg > 0 $$
(for $  \abg < 0 $ the transformation matrix would be complex, not real,
 for $  \abg = 0$ it would be singular).
Therefore we have in the case $\abg >0 $ a one--parametric set of non--equivalent Manin triples 
$$  [\tilde{X}^1,\tilde{X}^2]=- b \tilde{X}^2, \, 
[\tilde{X}^2,\tilde{X}^3] = 0 , \; 
[\tilde{X}^3,\tilde{X}^1] = b \tilde{X}^3, \; b >0  $$
and we must find representants of remaining classes of possible Manin triples.
We choose the forms
$$ [\tilde{X}^1,\tilde{X}^2]=0, \, 
[\tilde{X}^2,\tilde{X}^3] = b  \tilde{X}^2 , \; 
[\tilde{X}^3,\tilde{X}^1] = - b \tilde{X}^1, \; b >0  $$
 for   $  \abg < 0 $  and 
$$ [\tilde{X}^1,\tilde{X}^2]= \tilde{X}^2, \, 
[\tilde{X}^2,\tilde{X}^3] =   \tilde{X}^2  , \; 
[\tilde{X}^3,\tilde{X}^1] = - ( \tilde{X}^1+\tilde{X}^3 ) $$
for  $  \abg = 0, \, \alpha \neq 0 \vee \beta \neq 0 \vee \gamma \neq 0 $
and easily verify that every possible 2nd subalgebra $\tcg$ can be taken to one 
of the given forms by transformation (\ref{tfnb}) which doesn't change the structure constants 
of the 1st subalgebra $\cg=sl(2,\Re)$.\smallskip

Details of computations for each Bianchi algebra are given in the Appendix.

\section{Results: 6--dimensional  Manin triples}

The forms of non--equivalent Manin triples were choosen according to the following criteria:
The 1st subalgebra is in the Bianchi form, the 2nd is in the form closest to Bianchi, i.e.
Bianchi form if possible, or the structure constants are multiple of the Bianchi ones, or 
form a permution of the Bianchi ones, or, if neither is possible, are choosen to
be as many zeros and small integers as possible.

In order to shorten the list, we have not explicitly written out the structure of algebras 
that can be found by the duality transform $\cg \leftrightarrow \tcg $ from the 
ones given in the list.

\begin{enumerate}
\item Dual algebras to Bianchi algebra  $IX$:
$$  [X_1,X_2]=X_3, \; [X_2,X_3] = X_1, \; [X_3,X_1] = X_2. $$
Dual algebras:
\begin{enumerate} 
\item Bianchi algebra $I$
$$ [\tilde{X}^1,\tilde{X}^2]= 0 ,  \, 
[\tilde{X}^2,\tilde{X}^3] = 0 , \; 
[\tilde{X}^3,\tilde{X}^1] = 0. $$ 
\item Bianchi algebra $V$
$$ [\tilde{X}^1,\tilde{X}^2]=- b \tilde{X}^2, \, 
[\tilde{X}^2,\tilde{X}^3] = 0 , \; 
[\tilde{X}^3,\tilde{X}^1] = b \tilde{X}^3, \; b > 0  . $$ 
\end{enumerate}

\item Dual algebras to Bianchi algebra  $VIII$:
$$  [X_1,X_2]=-X_3, \; [X_2,X_3] = X_1, \; [X_3,X_1] = X_2. $$
Dual algebras: 
\begin{enumerate} 
\item Bianchi algebra $I$
$$[\tilde{X}^1,\tilde{X}^2]= 0 ,  \, 
[\tilde{X}^2,\tilde{X}^3] = 0 , \; 
[\tilde{X}^3,\tilde{X}^1] = 0. $$ 
\item Bianchi algebra $V$
\begin{enumerate}
\item $ [\tilde{X}^1,\tilde{X}^2]=- b \tilde{X}^2, \, 
[\tilde{X}^2,\tilde{X}^3] = 0 , \; 
[\tilde{X}^3,\tilde{X}^1] = b \tilde{X}^3, \; b >0 . $ 
\item $ [\tilde{X}^1,\tilde{X}^2]=0, \, 
[\tilde{X}^2,\tilde{X}^3] = b  \tilde{X}^2 , \; 
[\tilde{X}^3,\tilde{X}^1] = - b \tilde{X}^1, \; b >0 . $ 
\item $ [\tilde{X}^1,\tilde{X}^2]= \tilde{X}^2, \, 
[\tilde{X}^2,\tilde{X}^3] =   \tilde{X}^2  , \; 
[\tilde{X}^3,\tilde{X}^1] = - ( \tilde{X}^1+\tilde{X}^3 ) . $
\end{enumerate}
\end{enumerate}

\item Dual algebras to Bianchi algebra  $VII_{a}$:
$$  [X_1,X_2]=-a X_2+X_3, \; [X_2,X_3] = 0, \; [X_3,X_1] = X_2+ a X_3, 
\; a >0. $$
Dual algebras:
\begin{enumerate} 
\item Bianchi algebra $I$
$$ [\tilde{X}^1,\tilde{X}^2]= 0, \, 
[\tilde{X}^2,\tilde{X}^3] = 0 , \; 
[\tilde{X}^3,\tilde{X}^1] = 0. $$ 
\item Bianchi algebra $II$
\begin{enumerate}  
\item $ [\tilde{X}^1,\tilde{X}^2]=0, \, 
[\tilde{X}^2,\tilde{X}^3] = \tilde{X}^1 , \; 
[\tilde{X}^3,\tilde{X}^1] = 0 .$ 
\item $ [\tilde{X}^1,\tilde{X}^2]=0, \, 
[\tilde{X}^2,\tilde{X}^3] = - \tilde{X}^1 , \; 
[\tilde{X}^3,\tilde{X}^1] = 0 . $
\end{enumerate} 
\item Bianchi algebra $VII_{1/a}$
$$ [\tilde{X}^1,\tilde{X}^2]= b ( -  \frac{1}{a} \tilde{X}^2 + \tilde{X}^3), \, 
[\tilde{X}^2,\tilde{X}^3] = 0, \; 
[\tilde{X}^3,\tilde{X}^1] = b (\tilde{X}^2+ \frac{1}{a} \tilde{X}^3), 
\; b \in { \Re } - \{ 0 \} . $$ 
\end{enumerate}

\item Dual algebras to Bianchi algebra  $VII_{0}$:
$$  [X_1,X_2]=0, \; [X_2,X_3] = X_1, \; [X_3,X_1] = X_2. 
 $$
Dual algebras: 
\begin{enumerate} 
\item Bianchi algebra $I$
$$ [\tilde{X}^1,\tilde{X}^2]= 0, \, 
[\tilde{X}^2,\tilde{X}^3] = 0 , \; 
[\tilde{X}^3,\tilde{X}^1] = 0. $$ 
\item Bianchi algebra $II$
\begin{enumerate}  
\item $ [\tilde{X}^1,\tilde{X}^2]= \tilde{X}^3, \, 
[\tilde{X}^2,\tilde{X}^3] = 0 , \; 
[\tilde{X}^3,\tilde{X}^1] = 0 .$ 
\item $ [\tilde{X}^1,\tilde{X}^2]= - \tilde{X}^3, \, 
[\tilde{X}^2,\tilde{X}^3] = 0 , \; 
[\tilde{X}^3,\tilde{X}^1] = 0 . $
\end{enumerate} 
\item Bianchi algebra $IV$
$$ [\tilde{X}^1,\tilde{X}^2]= b ( - \tilde{X}^2 + \tilde{X}^3), \, 
[\tilde{X}^2,\tilde{X}^3] = 0, \; 
[\tilde{X}^3,\tilde{X}^1] = b  \tilde{X}^3, 
\; b \in { \Re } - \{ 0 \} . $$ 
\item Bianchi algebra $V$
\begin{enumerate}
\item $ [\tilde{X}^1,\tilde{X}^2]=  - \tilde{X}^2 , \, 
[\tilde{X}^2,\tilde{X}^3] = 0, \; 
[\tilde{X}^3,\tilde{X}^1] =  \tilde{X}^3, 
. $ 
\item $ [\tilde{X}^1,\tilde{X}^2]=  0 , \, 
[\tilde{X}^2,\tilde{X}^3] = b \tilde{X}^2, \; 
[\tilde{X}^3,\tilde{X}^1] = -b \tilde{X}^1, 
\, b >0  . $
\end{enumerate}

\end{enumerate}

\item Dual algebras to Bianchi algebra  $VI_{a}$:
$$  [X_1,X_2]=-a X_2-X_3, \; [X_2,X_3] = 0, \; [X_3,X_1] = X_2+ a X_3, 
\; a >0, \, a \neq 1 . $$
Dual algebras: 
\begin{enumerate} 
\item Bianchi algebra $I$
$$ [\tilde{X}^1,\tilde{X}^2]= 0, \, 
[\tilde{X}^2,\tilde{X}^3] = 0 , \; 
[\tilde{X}^3,\tilde{X}^1] = 0. $$ 
\item Bianchi algebra $II$
$$ [\tilde{X}^1,\tilde{X}^2]=0, \, 
[\tilde{X}^2,\tilde{X}^3] = \tilde{X}^1 , \; 
[\tilde{X}^3,\tilde{X}^1] = 0 . $$ 
\item Bianchi algebra $VI_{1/a}$
\begin{enumerate}
\item $ [\tilde{X}^1,\tilde{X}^2]=- b ( \frac{1}{a} \tilde{X}^2+\tilde{X}^3), \, 
[\tilde{X}^2,\tilde{X}^3] = 0, \; 
[\tilde{X}^3,\tilde{X}^1] = b (\tilde{X}^2+ \frac{1}{a} \tilde{X}^3), \; b \in { \Re } - \{ 0 \} . $ 
\item $ [\tilde{X}^1,\tilde{X}^2]= \tilde{X}^1 , \, 
[\tilde{X}^2,\tilde{X}^3] = \frac{a+1}{a-1} (\tilde{X}^2+\tilde{X}^3) , \; 
[\tilde{X}^3,\tilde{X}^1] = \tilde{X}^1 . $
\item $ [\tilde{X}^1,\tilde{X}^2]= \tilde{X}^1 , \, 
[\tilde{X}^2,\tilde{X}^3] = \frac{a-1}{a+1} (-\tilde{X}^2+\tilde{X}^3) , \; 
[\tilde{X}^3,\tilde{X}^1] = -\tilde{X}^1 .$ 
\end{enumerate}

\end{enumerate}

\item Dual algebras to Bianchi algebra  $VI_{0}$:
$$  [X_1,X_2]=0, \; [X_2,X_3] = X_1, \; [X_3,X_1] = -  X_2. 
 $$
Dual algebras: 
\begin{enumerate} 
\item Bianchi algebra $I$
$$ [\tilde{X}^1,\tilde{X}^2]= 0, \, 
[\tilde{X}^2,\tilde{X}^3] = 0 , \; 
[\tilde{X}^3,\tilde{X}^1] = 0. $$ 
\item Bianchi algebra $II$
 $$ [\tilde{X}^1,\tilde{X}^2]= \tilde{X}^3, \, 
[\tilde{X}^2,\tilde{X}^3] = 0 , \; 
[\tilde{X}^3,\tilde{X}^1] = 0 .$$ 
\item Bianchi algebra $IV$
\begin{enumerate}
\item $ [\tilde{X}^1,\tilde{X}^2]= b ( - \tilde{X}^2 + \tilde{X}^3), \, 
[\tilde{X}^2,\tilde{X}^3] = 0, \; 
[\tilde{X}^3,\tilde{X}^1] = b  \tilde{X}^3, 
\; b \in { \Re } - \{ 0 \} . $
\item $ [\tilde{X}^1,\tilde{X}^2]=  ( - \tilde{X}^1 + \tilde{X}^2 + \tilde{X}^3), \, 
[\tilde{X}^2,\tilde{X}^3] = \tilde{X}^3, \; 
[\tilde{X}^3,\tilde{X}^1] =  - \tilde{X}^3 
 . $
\end{enumerate} 
\item Bianchi algebra $V$
\begin{enumerate}
\item $ [\tilde{X}^1,\tilde{X}^2]=  - \tilde{X}^2 , \, 
[\tilde{X}^2,\tilde{X}^3] = 0, \; 
[\tilde{X}^3,\tilde{X}^1] =  \tilde{X}^3. $ 
\item $ [\tilde{X}^1,\tilde{X}^2]=  - \tilde{X}^1+ \tilde{X}^2 , \, 
[\tilde{X}^2,\tilde{X}^3] =  \tilde{X}^3, \; 
[\tilde{X}^3,\tilde{X}^1] = - \tilde{X}^3. $ 
\item $ [\tilde{X}^1,\tilde{X}^2]=  0 , \, 
[\tilde{X}^2,\tilde{X}^3] = - \tilde{X}^2, \; 
[\tilde{X}^3,\tilde{X}^1] =  \tilde{X}^1 . $ 
\end{enumerate}

\end{enumerate}

\item Dual algebras to Bianchi algebra  $V$:
$$  [X_1,X_2]=-X_2, \; [X_2,X_3] = 0, \; [X_3,X_1] = X_3. 
 $$
Dual algebras: 
\begin{enumerate} 
\item Bianchi algebra $I$
$$ [\tilde{X}^1,\tilde{X}^2]= 0, \, 
[\tilde{X}^2,\tilde{X}^3] = 0 , \; 
[\tilde{X}^3,\tilde{X}^1] = 0. $$ 
\item Bianchi algebra $II$
\begin{enumerate}
\item  $ [\tilde{X}^1,\tilde{X}^2]= 0, \, 
[\tilde{X}^2,\tilde{X}^3] = \tilde{X}^1 , \; 
[\tilde{X}^3,\tilde{X}^1] =  0  .$
\item  $ [\tilde{X}^1,\tilde{X}^2]= \tilde{X}^3, \, 
[\tilde{X}^2,\tilde{X}^3] = 0 , \; 
[\tilde{X}^3,\tilde{X}^1] =  0  .$
\end{enumerate}
\end{enumerate}
  and dual algebras ($\cg \leftrightarrow \tcg $) to algebras given above for  $VI_0$, $VII_0 $, $VIII$, $IX$.   

\item Dual algebras to Bianchi algebra  $IV$:
$$  [X_1,X_2]=-X_2+X_3, \; [X_2,X_3] = 0, \; [X_3,X_1] = X_3. 
 $$
Dual algebras: 
\begin{enumerate} 
\item Bianchi algebra $I$
$$ [\tilde{X}^1,\tilde{X}^2]= 0, \, 
[\tilde{X}^2,\tilde{X}^3] = 0 , \; 
[\tilde{X}^3,\tilde{X}^1] = 0. $$ 
\item Bianchi algebra $II$
\begin{enumerate}
\item  $ [\tilde{X}^1,\tilde{X}^2]= 0, \, 
[\tilde{X}^2,\tilde{X}^3] = \tilde{X}^1 , \; 
[\tilde{X}^3,\tilde{X}^1] =  0  .$
\item  $ [\tilde{X}^1,\tilde{X}^2]= 0, \, 
[\tilde{X}^2,\tilde{X}^3] = - \tilde{X}^1 , \; 
[\tilde{X}^3,\tilde{X}^1] =  0  .$
\item  $ [\tilde{X}^1,\tilde{X}^2]= 0, \, 
[\tilde{X}^2,\tilde{X}^3] = 0 , \; 
[\tilde{X}^3,\tilde{X}^1] = b \tilde{X}^2, \, b \in \Re - \{ 0 \} .$ 
\end{enumerate}
\end{enumerate}
and dual algebras ($\cg \leftrightarrow \tcg $) to algebras given above for $VI_0$, $VII_0$.  

\item Dual algebras to Bianchi algebra  $III$:
$$  [X_1,X_2]=-X_2-X_3, \; [X_2,X_3] = 0, \; [X_3,X_1] = X_2+X_3. $$
Dual algebras: 
\begin{enumerate} 
\item Bianchi algebra $I$
$$ [\tilde{X}^1,\tilde{X}^2]= 0, \, 
[\tilde{X}^2,\tilde{X}^3] = 0 , \; 
[\tilde{X}^3,\tilde{X}^1] = 0. $$ 
\item Bianchi algebra $II$
$$ [\tilde{X}^1,\tilde{X}^2]=0, \, 
[\tilde{X}^2,\tilde{X}^3] = \tilde{X}^1 , \; 
[\tilde{X}^3,\tilde{X}^1] = 0 . $$ 
\item Bianchi algebra $III$
\begin{enumerate}
\item $ [\tilde{X}^1,\tilde{X}^2]=- b (\tilde{X}^2+\tilde{X}^3), \, 
[\tilde{X}^2,\tilde{X}^3] = 0, \; 
[\tilde{X}^3,\tilde{X}^1] = b (\tilde{X}^2+\tilde{X}^3), \; b \in { \Re } - \{ 0 \} . $ 
\item $ [\tilde{X}^1,\tilde{X}^2]= 0 , \, 
[\tilde{X}^2,\tilde{X}^3] = \tilde{X}^2+\tilde{X}^3 , \; 
[\tilde{X}^3,\tilde{X}^1] =  0 . $
\item $ [\tilde{X}^1,\tilde{X}^2]= \tilde{X}^1 , \, 
[\tilde{X}^2,\tilde{X}^3] = 0 , \; 
[\tilde{X}^3,\tilde{X}^1] = - \tilde{X}^1 . $
\end{enumerate}

\end{enumerate}

\item Dual algebras to Bianchi algebra  $II$:
$$  [X_1,X_2]=0, \; [X_2,X_3] = X_1, \; [X_3,X_1] = 0.
 $$
Dual algebras: 
\begin{enumerate} 
\item Bianchi algebra $I$
$$ [\tilde{X}^1,\tilde{X}^2]= 0, \, 
[\tilde{X}^2,\tilde{X}^3] = 0 , \; 
[\tilde{X}^3,\tilde{X}^1] = 0. $$ 
\item Bianchi algebra $II$ 
\begin{enumerate}
\item  $ [\tilde{X}^1,\tilde{X}^2]= \tilde{X}^3, \, 
[\tilde{X}^2,\tilde{X}^3] = 0 , \; 
[\tilde{X}^3,\tilde{X}^1] =  0  .$
\item  $ [\tilde{X}^1,\tilde{X}^2]= -\tilde{X}^3, \, 
[\tilde{X}^2,\tilde{X}^3] = 0 , \; 
[\tilde{X}^3,\tilde{X}^1] =  0  .$
\end{enumerate}
\end{enumerate}
 and dual algebras ($\cg \leftrightarrow \tcg $) to algebras given above  for 
$III$, $IV$, $VI_0$, $VI_a$,  $VII_0$, $VII_a$.    

\item Dual algebras to Bianchi algebra  $I$:
$$  [X_1,X_2]=0, \; [X_2,X_3] = 0, \; [X_3,X_1] = 0 . $$
Dual algebras: all Bianchi algebras (in their Bianchi forms)
\end{enumerate}

\section{Conclusions}

We have classified 6--dimensional 
Manin triples or, equivalently, 3--dimensional Lie bialgebras. 
In computations computer algebra systems Maple V and Mathematica 4 were 
used for solving the sets of 
linear and quadratic equations that follow from the Jacobi identities and similarity transformations. 
The results were calculated independently in both systems and afterwards were checked one against the other.
The complete list consists of 78 classes of  Manin triples 
(if one considers dual Lie bialgebras equivalent, then the count is 44).
An open problem that remains is detecting the Manin triples that belong to the same Drinfeld double or, 
in other words, the classification of the 6-dimensional Drinfeld doubles.

One of interesting results is the number of possible Lie bialgebra  
structures for the algebra $VIII$, i.e. $sl(2,{\Re})$. In this case there are
up to rescaling 3 non--equivalent Manin triples.
As mentioned in the Introduction, to every Manin triple correspond a pair of Poisson--Lie T--dual models. 
Therefore, there should exist 3 different pairs of non-abelian Poisson-Lie T-dual models for $sl(2,{\Re})$. 
Only one of them appeared in the literature so far \cite{alkltse:qpl}. There is a natural
question whether these models are equivalent (i.e. whether they correspond to the  decomposition of one Drinfeld double)
and if they lead after quantisation to the same quantum model.

\section*{Appendix: Most general form of $\tcg$ of Manin triple with given $\cg$}
In this Appendix we present our computations in some detail. For each Bianchi algebra we give 
solutions of the mixed Jacobi identities (\ref{mjid}), i.e. linear equations in $\tilde{f}$, the 
remaining non--trivial Jacobi identities in $\tcg$ (\ref{jidd}), i.e. in general quadratic equations in  
$\tilde{f}$ and their solutions, in general depending on several parameters $\alpha,\beta,\ldots$ 
Finally we specify values of parameters allowing transformation (\ref{tfnb}) of $\tcg$ into 
forms of $\tcg$ given in the list of non--isomorphic Manin triples. 

\begin{itemize}
\item $\cg=IX$ \\
The mixed Jacobi identities (\ref{mjid})  
imply
$$\tilde{f^{23}}_3 = -\tilde{f^{12}}_1,  
\tilde{f^{23}}_2 = \tilde{f^{13}}_1, \tilde{f^{13}}_3 = \tilde{f^{12}}_2, \tilde{f^{23}}_1 = 0, 
\tilde{f^{12}}_3 = 0, \tilde{f^{13}}_2 = 0.$$
The  Jacobi identities in $\tcg$ (\ref{jidd}) in this case don't impose any new condition.
The general form of $\tcg$ is therefore
$$ [\tilde{X}^{1},\tilde{X}^{2}]=\alpha \tilde{X}^{1}+ \beta \tilde{X}^{2}, 
[\tilde{X}^{2},\tilde{X}^{3}]= \gamma \tilde{X}^{2} -\alpha \tilde{X}^{3}, 
[\tilde{X}^{3},\tilde{X}^{1}] = -\gamma \tilde{X}^{1} -\beta \tilde{X}^{3}$$
$\tcg$ can be transformed into 
\begin{itemize}
\item Bianchi algebra $I$ in the standard form $IX$ (a) if $\alpha=\beta=\gamma=0$, 
\item Bianchi algebra $V$ in the rescaled standard form $IX$ (b) with $b=\sqrt{\alpha^2+\beta^2+\gamma^2} $ 
otherwise.  
\end{itemize}

\item $\cg=VIII$ \\
The mixed Jacobi identities (\ref{mjid})  
imply
$$  \tilde{f^{12}}_1 = -\tilde{f^{23}}_3,  \tilde{f^{13}}_1 = \tilde{f^{23}}_2, 
\tilde{f^{12}}_2 = \tilde{f^{13}}_3, \tilde{f^{23}}_1 = 0, \tilde{f^{13}}_2 = 0, \tilde{f^{12}}_3 = 0. $$
The  Jacobi identities in $\tcg$ (\ref{jidd}) in this case don't impose any new condition.
The general form of $\tcg$ is therefore
$$ [ \tilde{X}^{1}, \tilde{X}^{2} ]= -\alpha \tilde{X}^{1} + \beta  \tilde{X}^{2}, 
\, [ \tilde{X}^{2}, \tilde{X}^{3}]= \gamma \tilde{X}^{2} + \alpha \tilde{X}^{3} , 
\, [ \tilde{X}^{3}, \tilde{X}^{1}]= -\gamma \tilde{X}^{1} - \beta \tilde{X}^{3}.$$ 
$\tcg$ can be transformed into 
\begin{itemize}
\item  Bianchi algebra $I$ in the standard form $VIII$ (a) if $\alpha=\beta=\gamma=0$,
\item Bianchi algebra $V$
\begin{itemize}
\item in the rescaled standard form $VIII$ (b) i. with $b=\sqrt{\abg}$ if $  \abg > 0 $, 
\item in the form $VIII$ (b) ii. with $b=\sqrt{-(\abg)}$ if $  \abg < 0 $, 
\item in the form $VIII$ (b) iii. if $ \abg = 0$, and $\alpha \neq 0 \vee \beta \neq 0 \vee \gamma \neq 0 
$\footnote{In order to avoid abundant parentheses, logical conjuctions written in terms of symbols are 
considered with higher priority than that written by words and, or.}.
\end{itemize}
\end{itemize}

\item $\cg=VII_{a}$ \\
The mixed Jacobi identities (\ref{mjid})  
imply
$$ \tilde{f^{13}}_2  = a\tilde{f^{13}}_3  , \tilde{f^{12}}_3  = -a\tilde{f^{13}}_3  ,  
\tilde{f^{23}}_3  = - \frac{a^2 \tilde{f^{23}}_2 +a^2 \tilde{f^{13}}_1  -\tilde{f^{23}}_2 
+\tilde{f^{13}}_1 }{2a} ,$$ 
$$\tilde{f^{12}}_1  = -\frac{a^2 \tilde{f^{23}}_2 +a^2 \tilde{f^{13}}_1  +\tilde{f^{23}}_2 
-\tilde{f^{13}}_1 }{2a}, \tilde{f^{12}}_2  = \tilde{f^{13}}_3  . $$
The  Jacobi identities in $\tcg$ (\ref{jidd}) reduce to
$$ 4 a\tilde{f^{23}}_1  \tilde{f^{13}}_3  + (a \tilde{f^{23}}_2)^2+2 a^2 \tilde{f^{23}}_2  \tilde{f^{13}}_1  
+(\tilde{f^{23}}_2)^2-2 \tilde{f^{23}}_2  \tilde{f^{13}}_1 +(a \tilde{f^{13}}_1) ^2 +(\tilde{f^{13}}_1)^2 
=0  .$$
The solutions of this equation give the following general forms of $\tcg$:
\begin{enumerate}
\item {} 
\begin{eqnarray*}
 [\tilde{X}^{1}, \tilde{X}^{2}]  & = &  -\frac{1}{2a}(a^2 \alpha+\beta a^2+\alpha-\beta) \tilde{X}^{1} +
\gamma  \tilde{X}^{2} - \gamma a \tilde{X}^{3} ,  \\
{[ \tilde{X}^{2}, \tilde{X}^{3} ]} & = &  
-\frac{1}{4\gamma a}(a^2 \alpha^2+2 \alpha \beta a^2+\alpha^2-2 \alpha \beta+
\beta^2 a^2+\beta^2) \tilde{X}^{1} + \alpha\tilde{X}^{2} \\
& & -\frac{1}{2a}(a^2 \alpha+\beta a^2-
\alpha+\beta) \tilde{X}^{3} ,  \\
{[ \tilde{X}^{3}, \tilde{X}^{1} ]} & = &  -\beta\tilde{X}^{1} -\gamma a \tilde{X}^{2} - \gamma \tilde{X}^{3}. 
\end{eqnarray*}
$\tcg$ can be transformed into 
\begin{itemize}
\item  Bianchi algebra $VII_{1/a}$ in the rescaled standard form $VII_{a}$ (c) with $b=-a \gamma$. 
\end{itemize}
\item
$ [ \tilde{X}^{1}, \tilde{X}^{2} ] = 0,
\, [ \tilde{X}^{2}, \tilde{X}^{3}] = \alpha \tilde{X}^{1} , 
\, [ \tilde{X}^{3}, \tilde{X}^{1}] = 0.$ \\
$\tcg$ can be transformed into 
\begin{itemize}
\item Bianchi algebra $I$ in the standard form $VII_{a}$ (a) if $ \alpha = 0$,
\item Bianchi algebra $II$
\begin{itemize}
\item in the standard form $VII_{a}$ (b) i. if $  \alpha > 0 $, 
\item in the form $VII_{a}$ (b) ii. if $  \alpha < 0 $. 
\end{itemize}
\end{itemize}
\end{enumerate}

\item $\cg=VII_{0}$
\\
The mixed Jacobi identities (\ref{mjid})  
imply
$$  \tilde{f^{12}}_1  = -\tilde{f^{23}}_3 , \tilde{f^{12}}_2  = \tilde{f^{13}}_3 , 
\tilde{f^{23}}_2  = \tilde{f^{13}}_1,\tilde{f^{13}}_2  = 0, \tilde{f^{23}}_1  = 0  . $$
The  Jacobi identities in $\tcg$ (\ref{jidd}) reduce to
$$ \tilde{f^{12}}_3 \tilde{f^{13}}_1 = 0.$$
The solutions of this equation give the following most general forms of $\tcg$:
\begin{enumerate}
\item {} 
\begin{eqnarray*}
[\tilde{X}^{1}, \tilde{X}^{2}]  & = & -\alpha  \tilde{X}^{1} +
\beta \tilde{X}^{2} + \gamma  \tilde{X}^{3} ,  \\
{[ \tilde{X}^{2}, \tilde{X}^{3} ]} & = &  \alpha  \tilde{X}^{3}, \\
{[ \tilde{X}^{3}, \tilde{X}^{1} ]} & = &  -\beta \tilde{X}^{3}. 
\end{eqnarray*}
$\tcg$ can be transformed into 
\begin{itemize}
\item  Bianchi algebra $I$ in the standard form $VII_{0}$ (a) if $\gamma  = \beta = \alpha  = 0$,
\item  Bianchi algebra $II$ 
\begin{itemize}
\item in the  form $VII_{0}$ (b) i.  
if $\gamma  >0 $ and $\beta =  \alpha  = 0$,
\item  in the  form $VII_{0}$ (b) ii.  
if $\gamma  < 0 $ and $\beta =  \alpha  = 0$,
\end{itemize}
\item  Bianchi algebra $IV$ in the rescaled standard form $VII_{0}$ (c) with 
$ b= - \frac{\beta^2+\alpha ^2 }{ \gamma} $ if $\gamma  \neq 0$ and $\beta \neq 0 \vee \alpha  \neq 0$,
\item  Bianchi algebra $V$ in the standard form $VII_{0}$ (d) i. with 
if $\gamma =0$ and $\beta \neq 0 \vee \alpha  \neq 0$.
\end{itemize}
\item
\begin{eqnarray*}
[\tilde{X}^{1}, \tilde{X}^{2}]  & = &  -\alpha  \tilde{X}^{1} +
  \beta  \tilde{X}^{2}  ,  \\
{[ \tilde{X}^{2}, \tilde{X}^{3} ]} & = &  
 \gamma  \tilde{X}^{2}  +  \alpha  \tilde{X}^{3} ,  \\
{[ \tilde{X}^{3}, \tilde{X}^{1} ]} & = &  -\gamma  \tilde{X}^{1}  - \beta  \tilde{X}^{3}. 
\end{eqnarray*}
$\tcg$ can be transformed into 
\begin{itemize}
\item Bianchi algebra $I$ in the standard form $VII_{0}$ (a) if $\alpha  = \beta  =  \gamma  = 0$,
\item Bianchi algebra $V$ 
\begin{itemize}
\item in the standard form $VII_{0}$ (d) i.  if $\gamma  = 0$,
\item in the form $VII_{0}$ (d) ii. with $b=|\gamma| $ if $\gamma  \neq 0$.
\end{itemize}
\end{itemize}
\end{enumerate}

\item $\cg=VI_{a}$
\\
The mixed Jacobi identities (\ref{mjid})  
imply
$$ \tilde{f^{13}}_1  = -\frac{-a^2 \tilde{f^{12}}_1 +a^2 \tilde{f^{23}}_3  -\tilde{f^{23}}_3 
-\tilde{f^{12}}_1 }{2a} , 
\tilde{f^{12}}_3  = a \tilde{f^{12}}_2 ,  \tilde{f^{13}}_2  = a \tilde{f^{12}}_2 ,$$ 
$$ \tilde{f^{13}}_3  = \tilde{f^{12}}_2 , \tilde{f^{23}}_2  = \frac{-a^2 \tilde{f^{12}}_1 + 
a^2 \tilde{f^{23}}_3  +\tilde{f^{23}}_3 +\tilde{f^{12}}_1 }{2a}   . $$
The  Jacobi identities in $\tcg$ (\ref{jidd}) reduce to
$$ 4 a \tilde{f^{23}}_1  \tilde{f^{12}}_2  +(a\tilde{f^{12}}_1)^2-2 a^2 \tilde{f^{12}}_1  \tilde{f^{23}}_3  
- 2 \tilde{f^{12}}_1  \tilde{f^{23}}_3 -(\tilde{f^{12}}_1)^2+(a \tilde{f^{23}}_3)^2 
-(\tilde{f^{23}}_3) ^2 =0 .$$
The solutions of this equation give the following most general forms of $\tcg$:
\begin{enumerate}
\item {} 
\begin{eqnarray*}
 {[\tilde{X}^{1}, \tilde{X}^{2}]}  & = & \alpha  \tilde{X}^{1}+ \beta \tilde{X}^{2}+ a \beta 
\tilde{X}^{3} ,  \\
{[ \tilde{X}^{2}, \tilde{X}^{3} ]} & = &  -\frac{a^2 \alpha ^2-2 \alpha  \gamma  a^2 -
 2 \alpha  \gamma -\alpha ^2+\gamma ^2 a^2-\gamma ^2}{4 a \beta } \tilde{X}^{1} \\
& & + 
\frac{(-a^2 \alpha +\gamma  a^2+\gamma +\alpha )}{2a} \tilde{X}^{2}+  \gamma  \tilde{X}^{3}
 ,  \\
{[ \tilde{X}^{3}, \tilde{X}^{1} ]} & = & 
\frac{-a^2 \alpha +\gamma  a^2-\gamma -\alpha }{2a} \tilde{X}^{1}
 - a \beta \tilde{X}^{2} -\beta \tilde{X}^{3} . 
\end{eqnarray*}
$\tcg$ can be transformed into 
\begin{itemize}
\item  Bianchi algebra $VI_{1/a}$ in the rescaled standard form $VI_{a}$ (c) with $b = -a\beta$. 
\end{itemize}
\item {}
\begin{eqnarray*}
{ [ \tilde{X}^{1}, \tilde{X}^{2} ]} & = & \alpha  \tilde{X}^{1},
\\ {[ \tilde{X}^{2}, \tilde{X}^{3}]} & = & \beta  \tilde{X}^{1} + 
\alpha \frac{a+1}{a-1} \tilde{X}^{2} + \alpha \frac{a+1}{a-1} \tilde{X}^{1}, 
\\ {[ \tilde{X}^{3}, \tilde{X}^{1}]} & = & \alpha  \tilde{X}^{1}.
\end{eqnarray*}
$\tcg$ can be transformed into 
\begin{itemize}
\item Bianchi algebra $I$ in the standard form $VI_{a}$ (a) if $\alpha = \beta  = 0$,
\item Bianchi algebra $II$ in the standard form $VI_{a}$ (b) if $\alpha = 0$ and $\beta  \neq 0$. 
\item Bianchi algebra $VI_{1/a}$ in the  form $VI_{a}$ (c) ii. if $\alpha  \neq 0$.
\end{itemize}
\item {}
\begin{eqnarray*}
  {[ \tilde{X}^{1}, \tilde{X}^{2} ]} & = & \alpha  \tilde{X}^{1},
\\ {[ \tilde{X}^{2}, \tilde{X}^{3}]} & = & \beta  \tilde{X}^{1}
 -\alpha  \frac{a-1}{a+1} \tilde{X}^{2} + \alpha  \frac{a-1}{a+1} \tilde{X}^{3}, 
\\ {[ \tilde{X}^{3}, \tilde{X}^{1}]} & = & -\alpha  \tilde{X}^{1}.
\end{eqnarray*}
$\tcg$ can be transformed into 
\begin{itemize}
\item Bianchi algebra $I$ in the standard form $VI_{a}$ (a) if $\alpha = \beta  = 0$,
\item Bianchi algebra $II$ in the standard form $VI_{a}$ (b) if $\alpha = 0$ and $\beta  \neq 0$. 
\item Bianchi algebra $VI_{1/a}$ in the form $VI_{a}$ (c) ii. if $\alpha  \neq 0$.
\end{itemize}
\end{enumerate}

\item $\cg=VI_{0}$\\
The mixed Jacobi identities (\ref{mjid})  
imply
$$ \tilde{f^{13}}_3  = \tilde{f^{12}}_2 , \tilde{f^{13}}_1  = \tilde{f^{23}}_2 , 
\tilde{f^{12}}_1  = -\tilde{f^{23}}_3 , \tilde{f^{13}}_2  = 0, \tilde{f^{23}}_1  = 0 . $$
The  Jacobi identities in $\tcg$ (\ref{jidd}) reduce to
$$ \tilde{f^{12}}_3  \tilde{f^{23}}_2 =0 .$$
The solutions of this equation give the following most general forms of $\tcg$:
\begin{enumerate}
\item {} 
\begin{eqnarray*}
[\tilde{X}^{1}, \tilde{X}^{2}]  & = & -\alpha  \tilde{X}^{1} + 
 \beta  \tilde{X}^{2}+ \gamma  \tilde{X}^{3} ,  \\
{[ \tilde{X}^{2}, \tilde{X}^{3} ]} & = &  \alpha  \tilde{X}^{3}  , \\
{[ \tilde{X}^{3}, \tilde{X}^{1} ]} & = &  -\beta  \tilde{X}^{3}. 
\end{eqnarray*}
$\tcg$ can be transformed into 
\begin{itemize}
\item  Bianchi algebra $I$ in the standard form $VI_{0}$ (a) if $\alpha = \beta  = \gamma  =0 $,
\item  Bianchi algebra $II$  in the  form $VI_{0}$ (b) if $ \gamma  \neq 0$ and $\alpha  = \beta  = 0$, 
\item  Bianchi algebra $IV$ 
\begin{itemize}
\item in the rescaled standard form $VI_{0}$ (c) i. with 
$b = \frac{\alpha ^2-\beta ^2}{\gamma} $
if $ \gamma  \neq 0$ and $\alpha ^2 \neq \beta ^2$,
\item in the form $VI_{0}$ (c) ii. with 
if $ \gamma  \neq 0$ and $\alpha ^2 = \beta ^2 \neq 0$,
\end{itemize}
\item  Bianchi algebra $V$ 
\begin{itemize}
\item in the standard form $VI_{0}$ (d) i.  if $\gamma =0$ and $\alpha ^2 \neq \beta ^2 $,
\item in the form $VI_{0}$ (d) ii.  if $\gamma =0$ and $\alpha ^2 = \beta ^2 \neq 0$,
\end{itemize}
\end{itemize}
\item {} 
\begin{eqnarray*}
{[\tilde{X}^{1}, \tilde{X}^{2}]}  & = & -\alpha  \tilde{X}^{1}+ \beta  \tilde{X}^{2}  ,  \\
{[ \tilde{X}^{2}, \tilde{X}^{3} ]} & = & \gamma  \tilde{X}^{2} + \alpha  \tilde{X}^{3}  , \\
{[ \tilde{X}^{3}, \tilde{X}^{1} ]} & = & -\gamma  \tilde{X}^{1} -\beta  \tilde{X}^{3} . 
\end{eqnarray*}
$\tcg$ can be transformed into 
\begin{itemize}
\item  Bianchi algebra $I$ in the standard form $VI_{0}$ (a) if $ \alpha  = \beta  =  \gamma  =0$,
\item  Bianchi algebra $V$ 
\begin{itemize}
\item in the  form $VI_{0}$ (d) i.  
if $\gamma  = 0$ and $\alpha ^2 \neq \beta ^2$,
\item in the  form $VI_{0}$ (d) ii.  
if $\gamma  = 0$ and $\alpha ^2 = \beta ^2$,
\item  in the  form $VI_{0}$ (d) iii. with $b = | \gamma | $ 
if $\gamma  \neq 0$.
\end{itemize}
\end{itemize}
\end{enumerate}

\item $\cg=V$ 
\\
The mixed Jacobi identities (\ref{mjid})  
imply
$$ \tilde{f^{12}}_1  = \tilde{f^{23}}_3 , \tilde{f^{13}}_3  = -\tilde{f^{12}}_2 , 
\tilde{f^{23}}_2  = -\tilde{f^{13}}_1 .$$
The  Jacobi identities in $\tcg$ (\ref{jidd}) in this case don't impose any new condition.
The general form of $\tcg$ is therefore
\begin{eqnarray*}
{[\tilde{X}^{1}, \tilde{X}^{2}]}  & = &  \alpha  \tilde{X}^{1} + \beta \tilde{X}^{2}+ 
 \gamma  \tilde{X}^{3},  \\
{[ \tilde{X}^{2}, \tilde{X}^{3} ]} & = &  \delta  \tilde{X}^{1} -\epsilon  \tilde{X}^{2}+ 
\alpha  \tilde{X}^{3}  , \\
{[ \tilde{X}^{3}, \tilde{X}^{1} ]} & = & -\epsilon  \tilde{X}^{1} -\zeta \tilde{X}^{2} + 
\beta \tilde{X}^{3}. 
\end{eqnarray*}
Finding the Bianchi forms of this algebra for all values of parameters seems to be rather complicated, 
because this case contains also all 2nd subalgebras of duals of Manin triples 
$(\cd,{\cg},{\tcg})$ with $\tcg \equiv V$ given above. 
Therefore we compute only the values of parameters for which $\tcg$ is isomorphic to 
$I, \ldots , V$. We find that $\tcg$ can be transformed into\footnote{It is helpful to exploit the fact 
that the commutant of $II$ is one--dimensional, i.e. suitably written matrix of structure coefficients 
has rank 1.}
\begin{itemize}
\item  Bianchi algebra $I$ in the standard form $V$ (a) if $ \alpha  = \beta = \gamma  = 
\delta  = \epsilon  = \zeta = 0 $,
\item  Bianchi algebra $II$  
\begin{itemize}
\item in the form $V$ (b) i. if 
\begin{itemize}
\item $\exists x,y$ s.t. $ \alpha  = x  \gamma , \beta = y  \gamma , 
\epsilon  = -x  y  \gamma , \zeta = - y^2 \gamma , \delta  = x^2 \gamma , \gamma  \neq 0 $
\item or $\alpha =\beta=\gamma =0$ and $\exists x$ s.t. $\epsilon  = -x \delta , 
\zeta = -x^2 \delta , x \neq 0, \delta  \neq 0 $ 
\item or $ \alpha  = \beta = \gamma  = 
\delta  = \epsilon  = 0, \zeta \neq 0$, 
\end{itemize}
\item in the  form $V$ (b) ii. if $ \alpha =\beta = \gamma =\epsilon  = \zeta = 0 $ 
 and $\delta  \neq 0$.
 \end{itemize}
\end{itemize}

\item $\cg=IV$ \\
The mixed Jacobi identities (\ref{mjid})  
imply
$$ \tilde{f^{12}}_3  = 0, \tilde{f^{12}}_2  = 0, \tilde{f^{23}}_2  = -\tilde{f^{13}}_1 -2 \tilde{f^{12}}_1 , 
\tilde{f^{23}}_3  = \tilde{f^{12}}_1 , \tilde{f^{13}}_3  = 0  . $$
The  Jacobi identities in $\tcg$ (\ref{jidd}) reduce to
$$ (\tilde{f^{12}}_1 )^2 = 0 .$$
The solution of this equation gives the most general form of $\tcg$:
\begin{eqnarray*}
{[\tilde{X}^{1}, \tilde{X}^{2}]}  & = &  0,  \\
{[ \tilde{X}^{2}, \tilde{X}^{3} ]} & = &  \alpha  \tilde{X}^{1} -\beta  \tilde{X}^{2}   , \\
{[ \tilde{X}^{3}, \tilde{X}^{1} ]} & = &  -\beta  \tilde{X}^{1} -\gamma  \tilde{X}^{2}. 
\end{eqnarray*}
$\tcg$ can be transformed into 
\begin{itemize}
\item  Bianchi algebra $I$ in the standard form $IV$ (a) if $\alpha = \beta  = \gamma  = 0$,
\item  Bianchi algebra $II$  
\begin{itemize}
\item in the standard form $IV$ (b) i. if $\gamma  = \beta =0 $ and  $\alpha >0$, 
\item in the  form $IV$ (b) ii. if $\gamma  = \beta = 0$ and $\alpha <0$, 
\item in the  form $IV$ (b) iii. with $b = - \gamma $ if $\gamma  \neq 0$ and 
$ \beta ^2 + \alpha  \gamma =0$, 
\end{itemize}
\item  Bianchi algebra $VI_{0}$ 
\begin{itemize}
\item in the rescaled standard form with 
$b = \gamma $ if $\gamma  \neq 0$ and $\beta ^2 + \alpha  \gamma  > 0 $. The corresponding Manin 
triple is dual to the triple $VI_{0}$ (c) i.
\item in the form 
$$
{[\tilde{X}^{1}, \tilde{X}^{2}]}   =   0,  \,
{[ \tilde{X}^{2}, \tilde{X}^{3} ]}  =    \tilde{X}^{2}, \,
{[ \tilde{X}^{3}, \tilde{X}^{1} ]}  =    \tilde{X}^{1}  $$
if $\gamma  = 0$ and $ \beta  \neq 0 $. The corresponding Manin 
triple is dual to the triple $VI_{0}$ (c) ii.
\end{itemize}
\item  Bianchi algebra $VII_{0}$ in the rescaled standard form with 
$b = \gamma $ if $\gamma  \neq 0$ and $\beta ^2 + \alpha  \gamma  < 0 $. The corresponding Manin 
triple is dual to the triple $VII_{0}$ (c) i.
\end{itemize}

\item $\cg=III$\\
The mixed Jacobi identities (\ref{mjid})  
imply
$$ \tilde{f^{13}}_3  = \tilde{f^{12}}_2 , \tilde{f^{12}}_1  = \tilde{f^{13}}_1 , 
\tilde{f^{12}}_3  = \tilde{f^{12}}_2 , \tilde{f^{13}}_2  = \tilde{f^{12}}_2 , 
\tilde{f^{23}}_3  = \tilde{f^{23}}_2  . $$
The  Jacobi identities in $\tcg$ (\ref{jidd}) reduce to
$$  \tilde{f^{23}}_1  \tilde{f^{12}}_2 - \tilde{f^{13}}_1  \tilde{f^{23}}_3  =0 .$$
The solutions of this equation give the following most general forms of $\tcg$:
\begin{enumerate}
\item {} 
\begin{eqnarray*}
[\tilde{X}^{1}, \tilde{X}^{2}]  & = &   \alpha  \tilde{X}^{1} +  \beta  \tilde{X}^{2} + 
\beta  \tilde{X}^{3} ,  \\
{[ \tilde{X}^{2}, \tilde{X}^{3} ]} & = & \frac{\alpha \gamma}{\beta}  \tilde{X}^{1} + \gamma \tilde{X}^{2}
+ \gamma \tilde{X}^{3} , \\
{[ \tilde{X}^{3}, \tilde{X}^{1} ]} & = & -\alpha  \tilde{X}^{1} -\beta  \tilde{X}^{2} 
- \beta  \tilde{X}^{3}. 
\end{eqnarray*}
$\tcg$ can be transformed into 
\begin{itemize}
\item  Bianchi algebra $III$ in the rescaled standard form $III$ (c) i. with 
$ b= 1/\beta  $.
\end{itemize}
\item {} 
\begin{eqnarray*}
[\tilde{X}^{1}, \tilde{X}^{2}]  & = & 0  ,  \\
{[ \tilde{X}^{2}, \tilde{X}^{3} ]} & = & \alpha  \tilde{X}^{1} + \beta  \tilde{X}^{2} + 
\beta  \tilde{X}^{3} , \\
{[ \tilde{X}^{3}, \tilde{X}^{1} ]} & = & 0. 
\end{eqnarray*}
$\tcg$ can be transformed into 
\begin{itemize}
\item  Bianchi algebra $I$ in the standard form $III$ (a) if $ \alpha  = \beta =0  $, 
\item  Bianchi algebra $II$ in the  form $III$ (b) i.  if $\beta  = 0 $ and $ \alpha  \neq 0 $,
\item  Bianchi algebra $III$ in the form $III$ (c) ii. if $\beta  \neq 0$.
\end{itemize}
\item {} 
\begin{eqnarray*}
[\tilde{X}^{1}, \tilde{X}^{2}]  & = & \alpha  \tilde{X}^{1} ,  \\
{[ \tilde{X}^{2}, \tilde{X}^{3} ]} & = & \beta \tilde{X}^{1} , \\
{[ \tilde{X}^{3}, \tilde{X}^{1} ]} & = & -\alpha  \tilde{X}^{1}. 
\end{eqnarray*}
$\tcg$ can be transformed into 
\begin{itemize}
\item  Bianchi algebra $I$ in the standard form $III$ (a) if $ \alpha  = \beta=0  $, 
\item  Bianchi algebra $II$ in the  form $III$ (b) i.  if $ \alpha  = 0 $ and $ \beta \neq 0 $,
\item  Bianchi algebra $III$ in the form $III$ (c) iii. if $ \alpha \neq 0$.
\end{itemize}
\end{enumerate}

\item $\cg=II$ \\
Finding the Bianchi forms of the 2nd algebra for all values of parameters again seems to 
be rather complicated, because it contains also all 2nd subalgebras of duals of Manin triples 
$(\cd,{\cg},\tcg)$ with $\tcg \equiv II$ given above. 
Therefore we compute only the values of parameters for which possible $\tcg$s are isomorphic to 
$I, II$. 

The mixed Jacobi identities (\ref{mjid})  
imply
$$  \tilde{f^{13}}_1  = \tilde{f^{23}}_2 , \tilde{f^{23}}_1  = 0, \tilde{f^{12}}_1  = -\tilde{f^{23}}_3 . $$
The  Jacobi identities in $\tcg$ (\ref{jidd}) reduce to
$$ -\tilde{f^{13}}_3  \tilde{f^{23}}_3 +\tilde{f^{23}}_3  \tilde{f^{12}}_2 -
2 \tilde{f^{12}}_3  \tilde{f^{23}}_2  =0, 
-2 \tilde{f^{13}}_2  \tilde{f^{23}}_3 -\tilde{f^{12}}_2  \tilde{f^{23}}_2 +\tilde{f^{23}}_2  
\tilde{f^{13}}_3 =0 .$$
The solutions of these equations give the following most general forms of $\tcg$:
\begin{enumerate}
\item {} 
\begin{eqnarray*}
{[\tilde{X}^{1}, \tilde{X}^{2}]}  & = &  -\alpha  \tilde{X}^{1} - \frac{2 \beta \alpha 
-\gamma  \delta }{
\gamma } \tilde{X}^{2} -\frac{\alpha ^2 \beta}{\gamma ^2} \tilde{X}^{3}
,  \\
{[ \tilde{X}^{2}, \tilde{X}^{3} ]} & = &   \gamma  \tilde{X}^{2}+ \alpha  \tilde{X}^{3}
, \\
{[ \tilde{X}^{3}, \tilde{X}^{1} ]} & = & -\gamma  \tilde{X}^{1} -\beta \tilde{X}^{2} 
- \delta  \tilde{X}^{3}
. 
\end{eqnarray*}
$\tcg$ of this form represents Bianchi algebras $IV,V$ only.

\item {} 
\begin{eqnarray*}
{[\tilde{X}^{1}, \tilde{X}^{2}]}  & = &   \alpha  \tilde{X}^{2} + \beta  \tilde{X}^{3}
,  \\
{[ \tilde{X}^{2}, \tilde{X}^{3} ]} & = &    0
, \\
{[ \tilde{X}^{3}, \tilde{X}^{1} ]} & = &  -\gamma  \tilde{X}^{2} -\delta  \tilde{X}^{3}
. 
\end{eqnarray*}
\begin{itemize}
\item  Bianchi algebra $I$ in the standard form $V$ (a) if $ \alpha = = \beta  = \gamma = \delta   =0$,
\item  Bianchi algebra $II$  
\begin{itemize}
\item in the form $II$ (b) i. if 
$\exists x: \, \gamma  = -x^2 \beta , \delta  = -x \beta , \alpha  = x \beta ,\beta  >0 $
 or $ \delta  =  \alpha  = \beta  =0, \gamma  < 0 $,
\item in the  form $II$ (b) ii. if
 $\exists x: \, \gamma  = -x^2 \beta , \delta  = -x \beta , 
\alpha  = x \beta ,\beta  <0 $ or $ \delta  =  \alpha  = \beta  =0, \gamma  > 0 $,
\end{itemize}
\end{itemize}
Bianchi algebras $III,IV,V,VI_{a},VI_0,VII_a,VII_0$ otherwise.

\item {} 
\begin{eqnarray*}
{[\tilde{X}^{1}, \tilde{X}^{2}]}  & = &  -\alpha  \tilde{X}^{1}+ \beta  \tilde{X}^{2}
+ \gamma  \tilde{X}^{1}
,  \\
{[ \tilde{X}^{2}, \tilde{X}^{3} ]} & = &    \alpha  \tilde{X}^{3}
, \\
{[ \tilde{X}^{3}, \tilde{X}^{1} ]} & = &   -\beta  \tilde{X}^{3}
. 
\end{eqnarray*}
\begin{itemize}
\item  Bianchi algebra $I$ in the standard form $V$ (a) if $ \alpha  = \beta  = \gamma  =0$,
\item  Bianchi algebra $II$  
\begin{itemize}
\item in the form $II$ (b) i. if $\alpha  = \beta  = 0, \gamma  >0$,
\item in the form $II$ (b) ii. if $\alpha  = \beta  = 0, \gamma  <0$,
\end{itemize}
\end{itemize}
Bianchi algebras $IV,V$ otherwise.
\end{enumerate}

\item $\cg=I$ \\
$\tcg $ might be any 3--dimensional Lie algebra, it can be brought to its 
Bianchi form by the transformation (\ref{tfnb}).

\end{itemize}


\begin{thebibliography}{1}

\bibitem{klse:dna}
C.Klim\v c\'{\i}k and P.\v Severa.
\newblock Dual non--{A}belian duality and the {D}rinfeld double.
\newblock {\em Phys.Lett. B}, 351:455--462, 1995.

\bibitem{kli:pltd}
C.Klim\v c\'{\i}k.
\newblock Poisson--{L}ie {T}-duality.
\newblock {\em Nucl.Phys B (Proc.Suppl.)}, 46:116--121, 1996.

\bibitem{hlasno:pltdm2dt}
L.~Hlavat\'y and L.~\v{S}nobl.
\newblock Poisson--{L}ie {T}--dual models with two--dimensional targets.
\newblock {\em e-print} hep-th/0110139.

\bibitem{iranci}
M.A. Jafarizadeh and A.Rezaei-Aghdam.
\newblock Poisson--{L}ie {T}-duality and {B}ianchi type algebras.
\newblock {\em Phys.Lett. B}, 458:470--490, 1999.

\bibitem{fof:sos}
J. M. Figueroa-O'Farrill. 
\newblock ${N}=2$ structures on solvable {L}ie algebras: the $c=9$ classification, 
\newblock {\em Comm. Math. Phys.}, 177:129, 1996.

\bibitem{gom:ctd}
X. Gomez.
\newblock Classification of three--dimensional Lie bialgebras.
\newblock {\em J. Math. Phys.}, 41:4939, 2000.

\bibitem{Drinfeld} V.G. Drinfeld. Quantum Groups. Proc. Int. Congr. Math. Berkeley, 798-820, 1986

\bibitem{Majid}
S. Majid.
\newblock Foundations of Quantum Group Theory.
\newblock Cambridge University Press, 2000

\bibitem{Landau}
L.D. Landau and E.M. Lifshitz. 
\newblock The Classical Theory of Fields.
\newblock Pergamon Press, 1987

\bibitem{alkltse:qpl}
A.Yu.Alekseev, C.Klim\v c\'{\i}k, and A.A.Tseytlin.
\newblock Quantum Poisson--{L}ie {T}-duality and {W}{Z}{N}{W} model.
\newblock {\em Phys.Lett. B}, 351:455--462, 1995.

\end{thebibliography}
\end{document}